\pdfoutput=1 
\documentclass[12pt]{article}
\usepackage{amsmath,amsfonts,amssymb,url,color,epsfig,graphics}
\usepackage[latin1]{inputenc}
\newcommand{\R}{{\mathbb{R}}}

\newcommand{\C}{{\mathbb{C}}}
\newcommand{\Z}{{\mathbb{Z}}}

\def\ha{\frac{1}{2}}

\def\pa{\partial}
\def\ra{\rightarrow}
\def\ga{\alpha}

\def\ge{\varepsilon}

\def\gg{\gamma}

\def\gl{\lambda}

\def\gs{\sigma}

\def\OPD{pseudo-differential operator}
\def\nghbd{neighbourhood }

\newtheorem{lemm}{Lemma}[section]
\newtheorem{prop}{Proposition}[section]

\newtheorem{theo}{Theorem}[section]

\begin{document}

\title{A proof of a  trace formula by Richard Melrose}

\author{Yves  Colin de Verdi\`ere\footnote{Universit\'e Grenoble-Alpes,
Institut Fourier,
 Unit{\'e} mixte
 de recherche CNRS-UGA 5582,
 BP 74, 38402-Saint Martin d'H\`eres Cedex (France);
{\color{blue} {\tt yves.colin-de-verdiere@univ-grenoble-alpes.fr}}}}

\maketitle
The goal of this note is to give  a new proof of the wave trace formula proved by Richard Melrose in the impressive paper \cite{Me-84}.
This trace formula is an extension of the Chazarain-Duistermaat-Guillemin
trace formula (denoted ``CDG trace formula'' in this paper)
to the case of a sub-Riemannian (``sR'') Laplacian on a 3D contact closed manifold. The  proof  uses a normal form constructed in the papers \cite{CHT-18,CHT-21},
following
the pioneering work \cite{Me-84}, 
in order to reduce to the case of the invariant Laplacian on the 3D-Heisenberg group.
We need also the propagation of singularities results of Victor Ivrii, Bernard Lascar and Richard Melrose \cite{Iv-76,La-82,Me-86}. 

% and the Toeplitz quantization of Louis Boutet de Monvel and Victor Guillemin \cite{B-G-81}.
{\it Acknowledgments: many thanks to Cyril for very useful comments!}

\section{The CDG  trace formula}
The following result was proved by Chazarain \cite{Ch-74} and refined by Duistermaat and Guillemin \cite{DG-75} (see also \cite{CdV-73,CdV-07} and Appendix
\ref{app:hist} for the history of the trace formulae): 
\begin{theo}
Let $(M,g)$ be a closed connected smooth  Riemannian manifold and $\gl_1=0 < \gl_2 \leq \cdots $ the spectrum of the Laplace operator, then we have  the following
equality of Schwartz distributions:
\[ \sum_{j=1}^\infty  e^{it\sqrt{\gl_j }}= T_0(t) + \sum_{\gg \in {\cal P}} T_\gg(t) ~{\rm mod}~ C^\infty \]
where ${\cal P}$ is the set of periodic geodesics, $ {\rm SingSupp}(T_0)=\{ 0 \} $ and, for $\gg $ a periodic geodesic, $ {\rm SingSupp}(T_\gg)\subset  \{ L_\gg \}
\cup \{- L_\gg \}$, where $L_\gg$ is the length of $\gg$.
Moreover, if $\gg $ is non degenerate,
\[ T_\gg (t)= \frac{L_0e^{i\pi  m(\gg)/2 }}{\pi {\rm det}({\rm Id}- P_\gg)^\ha}  (t+i0-L_\gg)^{-1}\left( 1+\sum_{j=1}^\infty  a_j (t-L_\gg)^j  \right) \]
where
\begin{itemize}
\item $L_0$ is the length of the primitive geodesic associated to $\gg$
\item $m(\gg)$ is the Morse index of $\gg$
\item $P_\gg$ is the linearized Poincaré map.
\end{itemize}
\end{theo}
This result gives a nicer proof of the main result of my thesis \cite{CdV-73} saying that, in the generic case, the Length spectrum, i.e. the set of lengths of closed geodesics,  is a spectral invariant.
This formula holds for any elliptic self-adjoint pseudo-differential operator $P$  of degree $1$ by replacing the periodic  geodesics by the
periodic  orbits of the Hamiltonian flow of the principal symbol of $P$ and  the Morse index  by a Maslov index.
Our goal is to prove that the same statement holds for sR Laplacians on a closed contact 3-manifold.
{\bf What we mean here  by `Melrose's trace formula'   is the
formulae of theorem 1, interpreted in this contact sR setting.}
Richard Melrose gave a proof which I found difficult (see Theorem 6.4 in \cite{Me-84}).

% \section{Toeplitz operators}

\section{Review of basic facts and notations}
For more details on this section, one can look at the paper \cite{CHT-18}.

\subsection{Contact 3D sR manifolds}

In what follows, $M$ is a closed (compact without boundary) connected manifold of dimension $3$ equipped with a smooth volume form $|dq|$.
We consider  also an oriented  contact distribution globally defined as the kernel of a non vanishing real valued 1-form $\ga$ so that 
$\ga \wedge d\ga $ is a volume form. Let  also $g$ be a metric  on the distribution $D=\ker \ga $. The ``co-metric'' $g^\star $ is defined
by $g^\star (q,p):= \| p_{|D_q} \|^2  $ where the norm is the dual norm of $g(q)$.
To such a set of data  is associated
\begin{itemize}
  \item 
    A geodesic flow denoted by $G_t,~t\in\R$: the Hamiltonian flow of $\sqrt{g^\star }$. The geodesics  are projections of the orbits of that flow onto $M$ and are everywhere tangent to $D$.
    We will often prefer to  consider the geodesic flow as the restriction to $g^\star =1$ of the Hamiltonian flow of $\ha g^\star$. 
  \item A Laplacian which is locally given by
    $\Delta = X_1^\star X_1+X_2^\star X_2 $ where $(X_1, X_2)$ is an orthonormal frame of $D$ and the adjoint is taken with respect to the measure $|dq|$.
  \item A canonical choice of a 1-form $\ga_g$ defining $D$ by assuming that  $d\ga_g$ restricted to $D$ is the oriented $g$-volume form on $D$.
  \end{itemize}

    From a famous theorem of H\"ormander, we know the  Laplacian is sub-elliptic and hence has a compact resolvent and a discrete spectrum
    $\gl_1=0 <\gl_2 \leq \cdots $ with smooth eigenfunctions. It follows from the sR  Weyl law
\[ \# \{ j | \gl_j \leq \gl \} \sim C \gl^2 \] 
that
    \[ {\rm Trace}(e^{it\sqrt{\Delta}}) = \sum _{j=1}^\infty e^{it\sqrt{\gl_j }} \]
    is a well defined Schwartz distribution often called the {\it wave trace} (see Appendix \ref{app:wave} for the link with the wave equation).
    Our goal is to extend the CDG formula to this case.

    The form $ \ga_g$ defines a {\it Reeb vector field}  $\vec{R} $ by the equations $\ga_g (\vec{R})= 1$ and $\iota (\vec{R})d\ga_g =0$.
    This vector field admits the following  Hamiltonian interpretation:
    the cone  $\Sigma =D^\perp $ is a symplectic sub-cone of $  T^\star M\setminus 0 $ \footnote{if $\Sigma =\{ (q,t\ga (q) |t\in\R \setminus 0, ~q\in M \}$,
     the symplectic form on $T^\star M $ restricts to $td\ga + dt\wedge \ga $ whose square is the volume form $-tdt\wedge \ga \wedge d\ga $}. 
We define the Hamiltonian
$\rho : \Sigma \ra \R$ by 
    $\rho (\ga )=\ga/\ga _g $ where $\ga \in \Sigma $ is a covector vanishing on $D$.
    The Hamiltonian vector field of $\rho $ is homogeneous of degree $0$ and the projection of this field onto  $M$ is  the Reeb vector field $\vec{R}$
    (see section 2.4 of \cite{CHT-18}). 

    \subsection{The 3D Heisenberg group $H_3$}  \label{sec:heis}

For this section, the reader could look at Section 3.1 of \cite{CHT-18}.
     We identify the Heisenberg group $H_3$ with $\R^3_{x,y,z}$ with the group law \[ (x,y,z)\star (x',y',z')=(x+x', y+y', z+z' +\ha(xy'-yx') )\] and the Lie algebra
is generated by $X,Y, Z$ with
\[ X =\pa_x +\ha y \pa _z,~Y=\pa _y -\ha x \pa_z,~Z=\pa _z \]
We choose $(D, g)$ by asking that $(X,Y)$ is an oriented orthonormal basis of $D$ for $g$ and $|dq|=|dxdy dz|$. 
We have $[X,Y]=-Z $ and the Reeb vector field is $Z$.
The Laplacian is $\Delta _3 =-(X^2+Y^2) $ and can be rewritten as
\begin{equation}\label{delta3}  \Delta _3= -|Z| \left( \left( \frac{X}{\sqrt{Z}} \right) ^2+ \left( \frac{Y}{\sqrt{Z}} \right)^2 \right) \end{equation}
on the complement of the kernel of  $Z$.  We write this as
\[ \Delta _3= |Z| \Omega \]
where $\Omega $ is an harmonic oscillator with spectrum $\{ 2l+1|l=0,1,\cdots \} $ (see \cite{CHT-18}, prop. 3.1).

We will need the following ``confining''  result (see \cite{Le-20}):
\begin{lemm} \label{lemm:confining}
Given $T>0$ , $\gs_0 \in \Sigma \setminus 0$ and $U$ a conic neighbourhood of $\gs_0$,  there exists a conic neighbourhood   $V \subset \subset U $ of $\gs_0$ so that $\forall t\in [-T,T]$, $G_t (V)\subset U $.
  \end{lemm}

\section{Speed of propagation}

First, we have the following
\begin{theo} If $u$ is a solution of an sR wave equation,
  \[\forall t\in\R,~  {\rm Support}(u(t))\subset B({\rm Support}(u(t=0))\cup {\rm Support}(du/dt (t=0)) , |t| )\]
  where $B(A,r)$ is the closed sR \nghbd~of radius $r$ of $A$.
\end{theo}
This result follows from the Riemannian case by passing to the limits (see Section 3 of \cite{Me-86}). 

We will also use the following  Theorem due to Victor Ivrii, Bernard Lascar and Richard Melrose \cite{Iv-76,La-82,Me-86} and revisited by Cyril Letrouit \cite{Le-21}):

\begin{theo} \label{theo:melrose} If $e(t,q,q')$ is the wave kernel of a sR Laplacian whose characteristic manifold $\Sigma $ is symplectic, i.e. $e$ is the Schwartz kernel of $\cos (t \sqrt{\Delta })$, then
  \[ WF'(e) \subset \{ (q,p,q',p',t,\tau )| \tau =\pm \sqrt{g^\star}, (q,p)=G_{\pm t} (q',p' ) \} \cup \{ (q,p,q,p,t, 0 )\} \]
  where $G_t$ is the geodesic flow. 
\end{theo}

\section{The local wave trace for the Heisenberg group}
As a preparation, we will prove the local trace formula for  the 3D-Heisenberg group $H_3$. 
The Laplacian $\Delta_3$  commutes with $Z$. We can hence use a partial Fourier decomposition of $L^2(H_3)$ identifying
it with the Hilbert integral 
\[ L^2(\R^3)= \int_\R^\oplus  {\cal H}_\zeta d\zeta \]
where ${\cal H}_\zeta :=\{ f|f(x,y,z+a)=e^{ia\zeta} f(x,y,z) \} $ which is identified to $L^2(\R^2)$ by looking at the value of $f$ at $z=0$.
In what follows, we will omit the space ${\cal H}_0$ which corresponds to a flat 2D-Euclidian Laplacian. 
Using this decomposition, the Laplacian rewrites as follows:
\[ \Delta_3 =\sum_{l=0} ^\infty (2l+1) \int_\R ^\oplus |\zeta | K_\zeta ^l  d\zeta \]
where the operator $K_\zeta ^l$ is the projector on the $l-$th Landau level of $\Delta _\zeta $, the restriction of $\Delta_3 $ to
${\cal H}_\zeta $, which is a magnetic Schr\"odinger operator
\[ H_\zeta := - \left( (\pa_x +i\zeta y/2)^2+(\pa_y -i\zeta x/2)^2 \right) \]
on $\R^2$ with magnetic field $\zeta dx\wedge dy $. 
The Schwartz kernel of  $K_\zeta ^l$ satisfies
\[K_\zeta^l(m,m)=\frac{|\zeta |}{2\pi }\]
(see Appendix \ref{app:K}).

Hence the half-wave operator writes 
\[ e^{it\sqrt{\Delta_3 }}= \sum _{l=0} ^\infty  \int_\R ^\oplus e^{it \sqrt{ (2l+1)|\zeta |} } K_\zeta  ^l d\zeta \]
  and  the local distributional  trace is given
  by
  \[ {\rm Trace }\left(e^{it\sqrt{\Delta_3 }} f \right) = \frac{1}{2\pi}\left( \int _{\R^3}f(q) |dq| \right)\sum _{l=0} ^\infty  \int_\R e^{it \sqrt{ (2l+1)|\zeta |} }|\zeta|
  d\zeta \]
  for $f\in C_0^\infty (H_3)$.
  This trace can be explicitly computed
   by  using the distributional Fourier transform 
  \[ \int _0^\infty e^{i\tau u} u^3 du = 6(\tau +i0)^{-4} .\]
    We get, for $t\ne 0$, 
  \[ {\rm Trace }\left(e^{it\sqrt{\Delta_3 }} f \right) = \frac{6}{\pi t^4 } \left( \sum _{l=0} ^\infty \frac{1}{(2l+1)^2} \right) \int _{H_3} f |dq|. \]
  In particular, this trace is smooth outside $t=0$ which is consistent with the fact that there is no periodic geodesic in $H_3$.

  The same result holds for the trace formula microlocalized near $\Sigma $:

  \begin{prop}\label{prop:microtrace}
    Let $P$ be a \OPD~ of degree $0$, which is compactly supported and so that
    $WF'(P)\cap \{ \zeta =0 \} =\emptyset $.
    Then
    \[ {\rm Trace}\left( e^{it\sqrt{\Delta_3 }} P \right) \]
      is smooth outside $t=0$. 
    \end{prop}
{\it Proof.--}
    We will first  prove  that
    %détails !! 
  \[ {\rm Trace}\left(e^{it\sqrt{\Delta_3 }}\chi \left(\frac{\Delta_3 }{|Z|^2},q\right)  \right)~, \]
  with $\chi \in C_0^\infty (\R\times H_3)$,  which rewrites  as
  \[ \frac{1}{2\pi } \sum _{l=0} ^\infty  \int _{H_3} |dq| \int_\R e^{it \sqrt{ (2l+1)|\zeta |} } \chi \left( \frac{2l+1}{|\zeta|} ,q\right) |\zeta| d\zeta ~, \]
  is smooth outside $t=0$.
  Let us define
  \[ I_l (t)=\int_\R e^{it \sqrt{ (2l+1)|\zeta |} } \chi \left( \frac{2l+1}{|\zeta|} ,q\right) |\zeta| d\zeta ~. \]
   We split the integral into  $\zeta >0$ and $\zeta <0$. Looking at the first integral and using the new variable
  $s=\sqrt{\zeta /2l+1}$, we get
  \[ I_l(t)= 2(2l+1)^2 \int_0^\infty  e^{it  (2l+1)s } \chi \left( \frac{1}{s^2},q\right)  s^3 ds  ~ \]
  Observing that the function $\chi \left( \frac{1}{s^2},q\right)  s^3$ is a smooth classical symbol supported away of $0$, integrations by part give
  that $I_l(t)=O\left(l^{-\infty }\right)$ as well as its derivatives. 
  
  Then we introduce
  \[ \bar{P} := \frac{1}{2\pi }\int_0^{2\pi } e^{it \Omega }P e^{-it \Omega } dt \]
  which is again a \OPD. We check that $P-\bar{P}=[ Q, \Omega ] $ for a  \OPD ~Q:
 if
 \[ S_t:=\int_0^{t } e^{is \Omega }(P-\bar{P}) e^{-is \Omega } dt, \]
 this holds true with
 \[ Q=\frac{1}{2\pi i}\int_0^{2\pi } e^{-it \Omega }S_t e^{it \Omega } dt~. \]  
The corresponding part of the trace vanishes,
  because $\Omega $ commutes with $\Delta_3 $. 
  Then the full symbol of $\bar{P}$ can be written as
  \[ \sum _{j=0}^\infty |\zeta |^{-j} p_j \left( \frac{I}{\zeta }, q \right) \] 
  with $p_j$ compactly supported.
  This  allows to reduce to the first case.
  \hfill $\square$

  \section{The trace formula for compact quotients of $H_3$}
  This section can be skipped. It contains  an example with a direct derivation of  Melrose's trace formula.
  
  Let us give a co-compact subgroup $\Gamma $ of $H_3$.  We will prove  Melrose's trace formula for $M:=\Gamma \backslash H_3$
  with Laplacian $\Delta_M$ which is $\Delta _3$, defined in Section \ref{sec:heis}, restricted to $\Gamma -$periodic functions.  We fix some time $T>0 $ and
  will look at the
  trace formula for $|t|\leq T $. We choose $\chi_D$ a smoothed fundamental domain, ie $\chi _D\in C_0^\infty (H_3)$ with
  $\sum _{\gamma \in \Gamma} \chi_D (\gg q)=1 $.
  We will denote by $\int_D\cdots $ the integral $\int_{H_3} \chi_D \cdots $. 

  We start with
  \[ e_M(t,q,q')=\sum _{\gg\in \Gamma } e_{3}(t,q,\gamma q' ) \]
  where $e_3$ is the half-wave kernel in $H_3$. 
  Because of the finite speed of propagation and the fact that $\Gamma $ is discrete, we have only to consider a finite sum for the trace:
  \[ {\rm Trace}\left(e^{it\sqrt{\Delta_M}}\right)= \int _D e_{3}(t,q,q) |dq| + \sum _{\gg \in \Gamma \setminus {\rm Id}, ~\min _q d(q,\gamma q)\leq T }\int _D e_{3}
  (t,q,\gamma q) |dq| \]
  The first term is smooth outside $t=0$ while the second one is given by the CDG trace formula.

With more details: let $c>0$ be small enough so that the cone $C_c=\{ g^\star < c \zeta^2 \}$ has the property $(\cal P)$ that there is no $\gamma \in
\Gamma \setminus {\rm Id}$ with a  geodesic from $q$ to $\gamma q$
of length  smaller than $T$ starting with Cauchy data in this cone.
  This is possible thanks to Lemma \ref{lemm:confining}.
  We can hence split the integrals
  $I_\gamma (t)=  \int _D e_{3} (t,q,\gamma q) |dq| $ into two pieces
  \[ I_\gg (t)={\rm Trace}\left( \left(e^{it \sqrt{\Delta_{3}} } \right) \tau _\gg \chi _D P \right) +
  {\rm Trace}\left(   \left(e^{it \sqrt{\Delta_{3}}}\right)\tau_\gg  \chi _D ({\rm Id}-P) \right) \]
      with $\tau _\gg (f) =f\circ \gg^{-1}$ and  $P= \psi (\Delta_3 /|Z|^2 ) $
      where  $\psi $ belongs to $ C_o^\infty (\R)$, is equal to $1$ near $0$ and is supported in  $ ]-{c},{c}[ $.
      The first term is smooth by  Theorem \ref{theo:melrose} and property $(\cal P)$. The second term corresponds to  the elliptic region
      and hence  we use the parametrix
      for the wave equation given by ``FIO''s
      as given in the CDG trace formula. We get then that the singularities of the wave trace locate on the   length spectrum.

  Note that the "heat trace" can be computed from the explicit expression of the spectrum. This is worked out in Appendix \ref{app:heat}. 

      \section{Normal forms}

      In what follows, $M$ is a closed 3D sR manifold of contact type equipped with a smooth volume. We denote by $\Delta $ the associated Laplacian. 
  The proof of the Melrose formula will be done by using a normal form allowing  a reduction to the case of Heisenberg.
  \subsection{Classical normal form}
  
  \begin{theo}
  Let $\Sigma $ be the characteristic manifold, i.e. the orthogonal of the distribution with respect to the duality,  and  
  let $\gs_0 \in \Sigma \setminus 0$, then there exists a conical \nghbd ~$U$  of $\gs_0$ and an homogeneous symplectic  diffeomorphism  $\chi $ of $U$
  onto a conical \nghbd
  of $(0,0, 0; 0,0, 1)$ in $T^\star H_3$, so that $g^\star_{H_3}\circ \chi = g^\star_M $.
  \end{theo}
  
  We use first   \cite{Me-84} (Prop. 2.3)  (or \cite{CHT-21} (Theorem 2.1)) to reduce to $\rho I $ where $\rho $ is the Reeb Hamiltonian
  and $I$ the harmonic oscillator Hamiltonian: this means that there is an homogeneous canononical transformation $\chi $ from a conic neighborhood of
  $\gs _0$ into $\Sigma_\sigma \times \R^2_{u,v}$ so that $\rho I\circ \chi = g^\star $ with $I=u^2+v^2$. 
  Let us denote by $\zeta$ the principal symbol of $Z$ in Equation (\ref{delta3}). 
  Then we use the normal form of Duistermaat-H\"ormander \cite{DH-72} (Prop. 6.1.3)  to reduce $\rho $ to $|\zeta |$ by a canonical transformation. 
  We get then the normal form $|\zeta | I $ which is the canonical decomposition of $g^\star_{H_3}$
  used in \cite{CHT-18} (see the principal symbols in Equation (\ref{delta3})).  
  \subsection{Quantum normal form}\label{sec:qnf}
  This is a 3-step reduction working in some conical \nghbd $C$ of a point of $\Sigma $. 
  \begin{enumerate}
   \item Using FIO's associated to $\chi $, we first reduce the Laplacian to
    a \OPD ~  of the form 
    $|Z| \Omega + R_0$ where  $R_0$
    is a \OPD ~  of degree $0 $.  This step is worked out in Theorem 5.2 of \cite{CHT-18}. 
  \item We can improve the previous normal form  so  that $R_0 $ commutes  with $\Omega $:  the cohomological equations
    \[ \{ |\zeta |I ,a \} = b \] (on $T^\star H_3$)
    where $b$,  vanishing on $\Sigma $ and  homogeneous of degree $j$ can be solved as shown in Appendix \ref{app:coh}  which is an improvement of what is proved 
    in \cite{CHT-18}.
        It follows that we get  a normal form
    $\Delta_3+ R_0$ with $R_0$ commuting with $\Omega $: the full symbol of $R_0$ is independent of  the angular part of the $(u,v)$ variables. 
  \item Using the spectral decomposition of $\Omega$, we get  a decomposition
    \[ \Delta \equiv  \oplus _{l=0}^\infty (2l+1) \Delta _l \Pi_l \]
    where the $\Delta _l$'s are \OPD s  of the form 
    \[ \Delta _l = \left(|Z| +\frac{1}{2l+1} R_0\right)   \]
    and $\Pi_l$ is the projector on the eigenspace of eigenvalue $2l+1$ of $\Omega$. 
    We can then use a reduction of the \OPD s  $|Z| +\frac{1}{2l+1} R_0$  to $|Z|$ by conjugating by elliptic \OPD s  ~ $A_l$
    depending smoothly of $\ge =1/(2l+1)$ and commuting with $\Omega $  as in \cite{DH-72}, proposition 6.1.4.
    We have
    \[ A_l^{-1} \left(|Z| +\frac{1}{2l+1} R_0 \right) A_l \equiv  |Z|  \]
    where $\equiv $ means modulo smoothing operators in $C$.  
  \end{enumerate}
  %Finally, we get:
  %\begin{theo}\label{theo:normal}
 %Using FIO's, the Laplace operator is microlocally near the characteristic manifold equivalent to the Heisenberg Laplacian.
   % More precisely, for any $\gs_0 \in  \Sigma $, there exists a conical \nghbd ~ $U$ of $\gs_0 $ and a Fourier integral operator $A$ microlocally  invertible in $U$, so
%that $\Delta = A^{-1}\Delta _3 A $ in $U$. 
  %  \end{theo} 

    \section{Proof of the Melrose trace formula}
Let us fix $T>0$ and try to prove Melrose's trace formula for times $t\in J:=[-T,T]$. 
Let us fix, for each $\sigma \in \Sigma $,  a conical \nghbd ~ $U_\gs $ of $\gs $ as in  section  \ref{sec:qnf}. 
We then take $W_\gs\Subset   V_\gs \Subset  U_\gs $ so that, for any $z\in V_\gs $ and any $t\in J $, 
$G_t (z) \in U_ \gs $ where $G_t$ is the geodesic flow. This is clearly possible using the classical normal form and the Lemma \ref{lemm:confining}. 
We then take a finite cover of $\Sigma $ by open cones $W_\ga:= W_{\gs_\ga} $ and a  finite pseudo-differential partition of unity
$(\chi_0, \chi_\ga (\ga  \in B)) $ so that $WF'(\chi_0)\cap \Sigma =\emptyset $, 
$\chi_\ga ={\rm Id } $ in $W_\ga $ and $WF'(\chi_\ga )\subset V_\ga $.
We have then to compute the traces
of $\left( \cos t \sqrt{\Delta } \right) \chi_0 $ and $\left( \cos t \sqrt{\Delta } \right) \chi_\ga $. 
We will prefer to use the wave equation now because the operator $\sqrt{\Delta }$ is not a \OPD !
We know from Theorem \ref{theo:melrose}  that, for $t\in J$, $WF'\left( \cos  (t\sqrt{\Delta} ) \chi_\ga \right) $ is a subset
of 
\[ \{ (z,z,t,0)|z\in V_\ga \}\cup \{  ( z,G_{\pm t}(z),t, \tau =\pm g^\star (z))|z\in V_\ga  \} .\] 
If $u(t) =\cos  (t\sqrt{\Delta})\chi_\ga  u_0 $, we have
$u_{tt} + \Delta u=0, u(0)=\chi_\ga u_0, u_t(0)=0.$
We can hence use the normal form and denote  by $\equiv $ the equality ``modulo smooth functions of $t\in J$'' to get
\[ Z_\ga (t):={ \rm Trace}(\cos (t\sqrt{\Delta } )\chi_\ga )\equiv { \rm Trace}(\cos (t\sqrt{\Delta_3+ R_0 } )\widetilde{\chi_\ga }),\]
where $ \widetilde{\chi_\ga }$ is the PDO  obtained by Egorov theorem when we take the normal form.
Then, because $R_0 $ commutes with $\Omega $, 
\[ Z_\ga (t)\equiv \sum _{l=0}^\infty { \rm Trace}\left(\cos \left(t\sqrt{(2l+1)\left(|Z| + \frac{1}{2l+1}R_0 \right)} \right)\Pi_l \widetilde{\chi_\ga} \right)\]
and
\[  Z_\ga (t)\equiv \sum _{l=0}^\infty { \rm Trace}\left( A_l^{-1}\cos \left(t\sqrt{(2l+1)|Z| } \right)A_l \Pi_l \widetilde{\chi_\ga} \right)\]
\[  Z_\ga (t)\equiv \sum _{l=0}^\infty { \rm Trace}\left(\cos \left(t\sqrt{(2l+1)|Z| } \right)A_l \Pi_l \widetilde{\chi_\ga} A_l^{-1}\right).\]
We can assume that $A_l$ is invertible on $WF'(\chi_\ga )$ and put $  \widetilde{\widetilde{\chi_\ga^l}}= A_l{\widetilde{\chi_\ga}}A_l^{-1}$. We get 
\[ Z_\ga (t)\equiv \sum _{l=0}^\infty { \rm Trace}\left(\cos \left(t\sqrt{(2l+1)|Z| } \right)A_l \Pi_l A_l^{-1}\widetilde{\widetilde{\chi_\ga^l}}\right)\]
Using the fact  $A_l$ commutes with $\Omega$ and hence with $\Pi_l$, we get finally
\[ Z_\ga (t)\equiv \sum _{l=0}^\infty { \rm Trace}\left(\cos \left(t\sqrt{(2l+1)|Z| } \right)\Pi_l\widetilde{\widetilde{\chi_\ga^l}} \Pi_l \right)\]
We can then apply a variant of the Proposition \ref{prop:microtrace}, more precisely of its proof,  where $P$ is replaced by $\oplus_{l=0}^\infty
\Pi_l \widetilde{\widetilde{\chi_\ga^l}}\Pi_l $ and using the fact that the $A_l$'s and hence the
$ \widetilde{\widetilde{\chi_\ga^l}}$ too
 are uniformly bounded \OPD s. 

It remains to study the part $Z_0(t)={\rm Trace}\left(\cos (t\sqrt{\Delta})\chi_0\right) $ which involves the elliptic part of the dynamics for which
we can use the FIO parametrix as in \cite{DG-75}. More precisely, for $t\in J $, the geodesic flow maps the microsupport of $\chi_0$ away of $\Sigma$, therefore
there exists a parametrix for $U(t)\chi_0 $ given by Fourier integral operators as in  the paper \cite{DG-75} and the calculation of the trace therefore follows the same
path. 
 \vfill \eject 

    \appendix
    {\bf \huge Appendices}
    
\section{Wave and half-wave equations} \label{app:wave}
Let $\Delta $ be  a self-adjoint positive sub-elliptic operator on a closed manifold. 
The wave equation is 
\[ \frac{\pa ^2u}{\pa t^2} + \Delta u=0, u(t=0)=u_0, \frac{\pa u}{\pa t} (t=0)=v_0 \]
This gives a one parameter group $U(t)=(U_0(t),U_1(t))$ on $L^2 \times L^2$. 
The trace of $U_0(t)$ is 
\[ Z_0(t)={\rm Trace}\left(\cos t\sqrt{\Delta } \right)=\sum _{j=1}^\infty \cos t\sqrt{\gl_j} \]
One can introduce also the half-wave equation $\frac{\pa u}{\pa t}=i\sqrt{\Delta }u,~ u(t=0)=u_0$.
The trace of the half-wave group is
$Z(t)=\sum_{j=1}^\infty  e^{it\sqrt{\gl_j}}$.
We have the relation 
$Z=2 H(Z_0)$ where $H$ is the $L^2$-projector multiplying the Fourier transform by the Heaviside function. 
It follows that the singularities of both distributions are easily related.

In the elliptic case, one can work directly with the half-wave group because $\sqrt{\Delta}$ is still an elliptic \OPD ~(Seeley's Theorem \cite{See-67}). 
This is no longer the case for sub-elliptic operators.

\section{The value of $K^l_\zeta (m,m)$} \label{app:K} 

Recall that $K_\zeta ^l $ is the orthogonal projector on the $l-$th Landau level  with a magnetic field in $\R^2$ equal to $\zeta dx\wedge dy$.
An easy rescaling shows that 
$K_\zeta ^l (m,m)=|\zeta | K_1 ^l (m,m) $.
We know from the Mehler formula (see \cite{Si-79}, p. 168)   that the heat kernel of the magnetic Schr\"odinger operator with constant magnetic field
equal to $1$ is given on the diagonal by
\[ e(t,m,m)=\frac{1}{4\pi \sinh t } \] 
On the other hand, we have
\[ e(t,m,m)=\sum _{l=0}^\infty e^{-(2l+1)t} K_1 ^l (m,m)\]
and
\[ \frac{1}{4\pi \sinh t } = \frac{1}{2\pi }\sum _{l=0}^\infty e^{-(2l+1)t} \]
Identifying both sums as Taylor series in $x=e^{-t}$ gives 
\[ K_1 ^l (m,m)= \frac{1}{2\pi } .\]

    \section{Toeplitz operators} \label{app:toeplitz}

    Let $\Sigma $ be a symplectic cone with a compact basis. Louis Boutet de Montvel and Victor Guillemin associate in  \cite{Bo-80,BG-81} to such a cone an
    Hilbert space and an algebra of operators called
    the Toeplitz operators with the same properties as the classical pseudo-differential operators. The latter case corresponds to the cone which is a cotangent cone.
    For an introduction, one can look at \cite{CdV-94}.
    
Two examples are implicitly present in this paper:

\begin{enumerate}
\item {\it Harmonic oscillator:} the harmonic oscillator $\Omega = -d_x^2 +x^2 $ is an elliptic self-adjoint Toeplitz operator; the cone $\Sigma $ is
$\R^2_{u,v}\setminus 0$ with the
symplectic form $du\wedge dv $ and the dilations
$\gl .(u,v)=(\sqrt{\gl }u, \sqrt{\gl}v) $. The symbol of $\Omega $ is $u^2+v^2$.

\item  {\it Quantization of the Reeb flow:} if $\Sigma \subset T^\star X \setminus 0$ is the characteristic cone of our sR Laplacian, one can quantize the
Reeb Hamiltonian $\rho $ as a first order elliptic Toeplitz operator of degree $1$. 
\end{enumerate}

    \section{A cohomological equation}\label{app:coh} 
    The following proposition is    a global formulation of the formal cohomological equations discussed in \cite{CHT-18} (section 5.1 and Appendix C) with a simple proof:
    \begin{prop} \label{prop:cohom} 
    We consider the cohomological equation
    \begin{equation} \label{equ:cohomo}
       \{ |\zeta | I ,A\} =B   \end{equation}
      where $A,B$ are smooth homogeneous functions in the cone $C:=\{ I < c |\zeta |\} $ with compact support in $q\in H_3$.
      If $B$ is homogeneous of degree $j$ and vanishes on $\Sigma :=\{ I=0 \} $,  Equation (\ref{equ:cohomo}) admits a solution $A$ homogeneous of degree $j-1$. 
    \end{prop}

    Restricting to $\zeta = 1$, reduces to prove the     
    following Lemma.
    \begin{lemm}
      Let us consider the differential equation
      \begin{equation} \label{equ:cohom}
        \frac{\pa a}{\pa \theta }+ \ha I  \frac{\pa a}{\pa z} = b(z,w ) \end{equation}
      with $(z,w)\in \R\times \{ |w|< c \}, w=|w|e^{i\theta }   $, $b$ smooth, compactly supported  in $z$ and $I=|w|^2$. We assume that $b(z,0)=0$.  
     
      Then Equation (\ref{equ:cohom}) admits a smooth solution $a$ depending smoothly of $b$. 
    \end{lemm}
    Any  smooth function $f$ in some disk in $\C$ admits a Fourier expansion 
    \[ f(w)=\sum_{n=0} ^\infty f_n(|w|^2) w^n + \sum_{n=1} ^\infty g_n(|w|^2) \bar{w}^n\]
    where the $f_n$'s and the $g_n$'s  are smooth \footnote{This is proved first in formal power series in $(w,\bar{w})$. This reduces the problem to functions
    which are flat at the origin. A Fourier expansion in $\theta $ allows to finish the proof.}.
We can use this expansion with  $f=a(z,.)$ and $g=b(z,.)$.  
We consider only the  sum of powers of $w^n$. The second can be worked out in a similar way. 
    We  then put $a=\sum _{n=1}^\infty a_n (z,I)w^n + a_0(z,I) $ and 
      $b=\sum _{n=1}^\infty b_n (z,I)w^n + I c_0(z,I) $. The factorization $b_0(z,I)= I c_0(z,I)$ follows from the assumption on $b$. 
    We can take   $a_0 (z,I)= 2\int_{-\infty }^z b_0(s,I) ds $. The equation for
    $a_n$,  with $n\geq 1$, writes
    \[ in a_n + \ha I \frac{\pa a_n }{\pa z}= b_n \]
    We can solve it, for $I\ne 0$,  by
    \[ a_n (z, I)= \frac{2}{I}\int_{-\infty}^0 b_n (z+u,I)e^{in 2u/I } du =\int_{-\infty}^0 b_n (z+Is/2,I)e^{in s } ds \]
    and $a_n (z,0)= b_n (z,0)/in $. We need to prove that $a_n$ is smooth: we check that $a_n$ is continuous and then the derivative are given by the same kind
    of integrals with derivatives of $b_n$ as shown by the last expression of $a_n$.
    The continuity of $a_n$ follows from an   integrations by parts in the first expression of $a_n$.

    Now we want to add  up the series  $\sum _n a_nw^n $. I was not able to do that directly and  will proceed as follows:
    the sum $\sum_n  a_n w^n $ is  convergent as a formal series along $\Sigma $,  because $a_n w^n =O(I^{n/2})$.
    Using Borel procedure, we need only to solve our cohomological equation with a flat righthandside.
    This follows clearly from the expression
    \[ a(z,w) =\int_{-\infty}^{0} b(z+It/2,e^{it}w) dt ~. \]
    In fact only the behaviour as $I\ra 0 $ could be a problem, but we can divide $b$ by any power of $I$. 
    
\section{The heat trace for a compact quotient of $H_3$}\label{app:heat}

Let us consider the discrete subgroup $\Gamma = \left( \sqrt{2\pi} \Z \right)^2 \times \pi \Z $ of $H_3$ identified with $\R^3_{x,y,z}$ as in Section \ref{sec:heis}. 
The spectrum of the sub-Laplacian defined in section \ref{sec:heis}  on $M=H_3/\Gamma $ is
the union of the spectrum of the flat torus $\R^2/ \left( \sqrt{2\pi} \Z \right)^2$ and
the eigenvalues $2m (2l+1), m\geq 1, ~l\geq 0,$ with multiplicities $2m$.
The corresponding part of the complexified heat trace is hence
\[ Z_o (z)= \sum _{m=1}^\infty 2m \sum _{l=0}^\infty e^{-2m (2l+1)z} \]
with  $\Re(z)>0$. Summing with respect to $l$ gives
\[  Z_o (z)= \sum _{m=1}^\infty \frac{m}{\sinh 2 m z } \]
which we rewrite as 
\[  Z_o (z)=\frac{1}{4z}\sum _{m\in \Z}^\infty \frac{2mz}{\sinh 2 m z }-\frac{1}{4z}\]
The Fourier transform of $\frac{x}{\sinh x} $ is
$\frac{\pi^2}{1+\cosh \pi \xi }$.
Applying Poisson summation formula to the last expression of $Z_o$, we get
\[ Z_o(z)= \frac{\pi^2}{16z^2}-\frac{1}{4z} + \frac{\pi^2}{4z^2} \sum _{n=1}^\infty \frac{1}{1 + \cosh \pi^2 n/ z } \] 
The first term gives the Weyl law. Each term in the sum w.r. to $n$ is equivalent to
\[ \frac{\pi^2}{2z^2} e^{-\pi^2 n/ z} \]
We observe that the lengths of the periodic geodesic of $M$ are the numbers $2\pi \sqrt{n} $. Hence we recover also the length spectrum giving contributions
of the order of
${\rm exp}(-L^2/4z )$ as in the Riemannian case as proved in \cite{CdV-73}.
It could be nice to derive an exact formula for the wave  trace from our expression of the heat trace. 
Similarly, the heat trace for the Riemannian Laplacian on $M$ was computed by Hubert Pesce \cite{Pe-94}.

%\section{Notations from \cite{CHT-18}}
    \section{A short history of the trace formulae} \label{app:hist}
    The trace formulae were first discovered  independently by two  groups of physicists: Martin Gutzwiller \cite{Gu-71} for a semi-classical Schr\"odinger operator
    and  Roger  Balian \& Claude Bloch in a
    very impressive series of papers
    for Laplacians in Euclidean domains \cite{Ba-Bl-70,Ba-Bl-71,Ba-Bl-72}.
    In \cite{Ba-Bl-72}, page 154, the authors
    suggested already a possible application to the  inverse spectral problems\footnote{They wrote ``The analysis of the eigenvalue density as a sum of oscillating
    terms gives a new insight into the problem
      of ``hearing the shape of a drum'' [Kac paper]. \dots It is convenient, for the discussion to start from the fact that the knowledge of eigenvalues determines
      uniquely the path generating function \dots
    Thus \dots the lengths of the closed stationary polygons are determined''. },
    an industry which just
    started at the end of the sixties. From the point of view of mathematics, the Poisson summation formula can be interpreted as a trace formula for the Euclidian
Laplacian on flat tori. Similarly, the famous Selberg trace formula \cite{Se-56}  (see also Heinz Huber \cite{Hu-59})  is a  trace formula for the Laplacian on
hyperbolic surfaces.
Then  my thesis \cite{CdV-73}, inspired by the work of Balian and Bloch and the Selberg trace formula, uses
the complex heat equation for general closed Riemannian manifold. The definitive  version, the CDG formula, using wave equation, was discovered by
Jacques Chazarain and the tandem  Hans Duistermaat \& Victor Guillemin in \cite{Ch-74,DG-75}. They use the power of 
 the Fourier Integral Operators calculus \cite{Ho-71,DH-72}. See \cite{CdV-07} for a review paper. Later results cover the cases of manifolds with boundaries and
 semi-classical versions. 

{\it Author references: no conflist of interest, no funding. }

\bibliographystyle{plain}

\end{document}